\documentstyle[12pt]{article}
    \newtheorem{th}{Theorem}
    
    \newtheorem{ex}{Example}
    \newtheorem{rem}{Remark}
    \newtheorem{dfn}{Definition}
    \newtheorem{prop}{Proposition}
    \newtheorem{cor}{Corollary}
    
    \newtheorem{lem}{Lemma}
    
    \newenvironment{prb}{\begin{quote}{\bf Problem.}}{\end{quote}}

\def\o{\otimes}

\def\vfi{\varphi}

\def\ss{\subset}

\def\bop{{} \oplus {}}

\def\hop{\hat{{} \oplus {}}}

\def\G{\Gamma}

\def\v{\vert}

    \newcommand{\operatorname}[1]{{\rm #1}\,}

\def\Ind{\operatorname{Ind}}

\def\Ch{\operatorname{Ch}}

\def\Tr{\operatorname{Tr}}

\def\Ker{\operatorname{Ker}}

\def\Im{\operatorname{Im}}

\def\Z{{\bf Z}}
\def\C{{\bf C}}

\def\vfi{\varphi}

\def\v{\vert}

\def\Ker{\operatorname{Ker}}

\def\Im{\operatorname{Im}}

\def\Z{{\bf Z}}
\newcommand{\Mat}[4]{\left( \begin{array}{cc}
                            #1 & #3 \\
                            #2 & #4
                      \end{array} \right)}
\newcommand{\Matt}[9]{\left( \begin{array}{ccc}
                            #1 & #2 & #3 \\
                            #4 & #5 & #6 \\
                            #7 & #8 & #9
                      \end{array} \right)}

\author{M.~Frank \and E.~Troitsky}
\title{Lefschetz Numbers and Geometry of Operators in W*-modules}

\begin{document}
\maketitle
\section{Introduction}

\noindent
The main goal of the present paper is to generalize the results
of~\cite{TroLNM,TroBoch} in the following way: To be able to define
$K_0(A)\o\C$-valued Lefschetz numbers of the first type of an endomorphism
$V$ on a C*-elliptic complex one usually assumes that $V=T_g$ for some
representation $T_g$ of a compact group $G$ on the C*-elliptic complex.
We try to refuse this restriction in the present paper. The price to pay
for this is twofold: \newline
(i) $\,$ We have to define Lefschetz numbers valued in some larger group
as $K_0(A)\o\C$.  \newline
(ii) We have to deal with W*-algebras instead of general unital
C*-algebras.

\noindent
To obtain these results we have got a number of by-product facts on
the theory of Hilbert W*- and C*-modules and on bounded module operators
on them which are of independent interest.

\smallskip \noindent
The present paper is organized as follows:
In \S 2 we prove the necessary facts on Hilbert W*-modules and
their bounded module mappings extending results of W.~L.~Paschke
\cite{Pa}, J.-F.~Havet \cite{Ha} and the first author
\cite{Frank:ZAA}.
In \S 3 we define Lefschetz numbers of two types and
show the main properties of them. In \S 4 we discuss the C*-case
and obstructions to refine the main results of \S 3.

\smallskip \noindent
Our standard references for the theory of Hilbert C*-modules
are the papers \cite{Pa,Pa2,DupFil,Lin:90/2,Frank:ZAA,Lin:92,Manu}
and the book of E.~C.~Lance \cite{Lance}. The topological
considerations are based on the publications
\cite{Mi,MF,TroGlobA,TroLNM,TroBoch,Manu}. We are going to
continue the investigations presented therein.

\newpage


\section{Hilbert W*-modules and module mappings}

\noindent
We want to show some more very nice properties of Hilbert W*-modules which
often do not appear in the general C*-case. This partial class of Hilbert
C*-modules was brought to the attention of the public by W.~L.~Paschke in his
classical paper \cite{Pa}, and they are of use in many cases. The facts below
can be reproved for the class of monotone complete C*-algebras carrying out
much technical work, cf.~\cite{Frank:MN}, but not for larger classes of
C*-algebras, in general. However, since we are going to understand the
structure of general Hilbert C*-modules and their C*-duals better it suffices
to
treat the W*-case, and we can avoid these technicalities. Let us start with
a property generalizing the (double) annihilator property of arbitrary subsets
of W*-algebras.

\begin{lem}  \label{lem1}
   Let $A$ be a W*-algebra and $\{ {\cal M}, \langle .,. \rangle \}$ be a
   Hilbert $A$-module. For every subset ${\cal S} \subseteq {\cal M}$ the
   bi-orthogonal set ${\cal S}^{\bot \bot} \subseteq {\cal M}$ is a Hilbert
   $A$-submodule and a direct summand of ${\cal M}$, as well as the orthogonal
   complement ${\cal S}^{\bot}$.
\end{lem}

\noindent
{\it Proof:\/}
The property of ${\cal S}^{\bot \bot} \subseteq {\cal M}$ to be a Hilbert
$A$-submodule is obvious by the definition of orthogonal complements.
Since the $A$-dual Banach $A$-module ${\cal M}'$ of ${\cal M}$ is a
self-dual Hilbert $A$-module by \cite[Th.~3.2]{Pa} one can consider the
Hilbert $A$-submodule ${\cal N}$ of ${\cal M}'$ consisting of the direct
sum of ${\cal S}^{\bot \bot} \hookrightarrow {\cal M}'$ and of the Hilbert
$A$-module of all $A$-linear bounded mappings from ${\cal M}$ to $A$
vanishing on ${\cal S}^{\bot \bot}$. The second summand is the orthogonal
complement of ${\cal S}^{\bot \bot}$ with respect to ${\cal M}'$ by
construction and hence, it is a self-dual Hilbert $A$-submodule and
direct summand of ${\cal N}$ by \cite[Th.~3.2,Th.~2.8]{Frank:ZAA}.
Consequently, the canonical embedding of ${\cal S}^{\bot \bot}$ into
${\cal N}$ is a direct summand of ${\cal N}$, and because of the
submodule inclusion ${\cal S}^{\bot \bot} \subseteq {\cal M}
\hookrightarrow {\cal N}$ it is a direct summand of ${\cal M}$, too.
$\: \bullet$

\smallskip \noindent
Example \ref{exam1} below shows that situations different to that described
at Lemma \ref{lem1} can appear e.~g.~for Hilbert C*-modules over the
C*-algebra $A=C([0,1])$.

\begin{lem}  \label{lem2}
    Let $A$ be a W*-algebra, $\{ {\cal M}, \langle .,. \rangle \}$ be a
    Hilbert $A$-module and $\phi$ be a bounded module operator on it. Then
    the kernel ${\rm Ker}(\phi)$ of $\phi$ is a direct summand of ${\cal M}$
    and has the property ${\rm Ker}(\phi) = {\rm Ker}(\phi)^{\bot\bot}$.
\end{lem}

\noindent
{\it Proof:\/} By \cite[Prop.~3.6]{Pa} every bounded module operator
$\phi$ on ${\cal M}$ continues to a bounded module operator on its
$A$-dual Hilbert $A$-module ${\cal M}'$. The kernel of the
extended operator is a direct summand of ${\cal M}'$ because of the
completeness of its unit ball with respect to the $\tau_2$-convergence
induced by the functionals $\{ f(\langle .,y \rangle) : f \in A_{*,1}, y
\in {\cal M}'\}$ there, (cf.~\cite[Th.~3.2]{Frank:ZAA}). Consequently,
the kernel of $\phi$ inside ${\cal M}$ has to coincide with its bi-orthogonal
complement in ${\cal M}$, and by Lemma \ref{lem1} it is a direct summand.
$\: \bullet$

\begin{ex}  \label{exam1} {\rm
   Note, that the kernel of bounded $A$-linear operators on Hilbert
   $A$-modules over arbitrary C*-algebras $A$ is not a direct summand, in
   general. For example, consider the C*-algebra $A=C([0,1])$ of all
   continuous functions on the interval [0,1] as a Hilbert $A$-module over
   itself equipped with the standard inner product $\langle a,b \rangle_A=
   ab^*$. Define the mapping $\phi_g$ by the formula $\phi_g(f)=g \cdot f$
   for a fixed function
   \[
   g(x)= \left\{ \begin{array}{ccc}
                 -2x+1 & : & x \leq 1/2 \\ 0 & : & x \geq 1/2
                 \end{array} \right.
   \]
   and for every $f \in A$. Then ${\rm Ker}(\phi_g)$ equals to the Hilbert
   $A$-submodule and (left) ideal $\{ f \in A : f(x)=0 \: \, {\rm for} \: \,
   x \in [0,1/2] \}$, being not a direct summand of $A$, but nevertheless,
   coinciding with the bi-orthogonal complement of itself with respect to
   $A$.       }
\end{ex}

\begin{cor} \label{corol1}
    Let $A$ be a W*-algebra, ${\cal M}$ and ${\cal N}$ be two Hilbert
    $A$-modules and $\phi : {\cal M} \rightarrow {\cal N}$ be a bounded
    $A$-linear mapping. Then the kernel ${\rm Ker}(\phi)$ of $\phi$ is a
    direct summand of ${\cal M}$ and has the property ${\rm Ker}(\phi) =
    {\rm Ker}(\phi)^{\bot\bot}$.
\end{cor}

\noindent
{\it Proof:\/} Consider the Hilbert $A$-module ${\cal K}$ formed as the
direct sum ${\cal K}={\cal M} \oplus {\cal N}$ equipped with the $A$-valued
inner product $\langle .,. \rangle_{{\cal M}}+\langle .,. \rangle_{{\cal N}}$.
The mapping $\phi$ can be identified with a bounded $A$-linear mapping
$\phi'$ on ${\cal K}$ acting on the direct summand ${\cal M}$ as $\phi$
and on the direct summand ${\cal N}$ as the zero operator. Since the
kernel of $\phi'$ is a direct summand of ${\cal K}$ containing ${\cal N}$
by Lemma \ref{lem2} its orthogonal complement is a direct summand of ${\cal
M}$.
The desired result turns out. $\: \bullet$

\medskip \noindent
Now we are in the position to give a description of the inner structure of
arbitrary Hilbert W*-modules generalizing an analogous statement for
self-dual Hilbert W*-modules by W.~L.~Paschke (\cite[Th.~3.12]{Pa}).

\begin{prop}  \label{prop1}
   Let $A$ be a W*-algebra and $\{ {\cal M}, \langle .,. \rangle \}$ be a
   left Hilbert $A$-module. Then ${\cal M}$ is the closure of a direct
   orthogonal sum of a family $\{ D_\alpha : \alpha \in I \}$ of norm-closed
   left ideals $D_\alpha \subseteq A$, where the closure of this direct sum
   is predetermined by the given on ${\cal M}$ $A$-valued inner product
   $\langle .,. \rangle$ and the $A$-valued inner products on the ideals are
   the standard $A$-valued inner product on $A$.
   Moreover, for every bounded $A$-linear mapping $r:{\cal M} \rightarrow A$
   there is a net $\{ x_\beta : \beta \in J \}$ of elements of ${\cal M}$
   for which the limit
   \[
   \|.\|_A-\lim_{\beta \in J} \langle y,x_\beta \rangle
   \]
   exists for every $y \in {\cal M}$ and equals $r(y)$.
\end{prop}

\noindent
{\it Proof:\/} Fix an arbitrary bounded $A$-linear mapping $r:{\cal M}
\rightarrow A$. The kernel of $r$ is a direct summand of ${\cal M}$ by
Corollary \ref{corol1}. Consider its orthogonal complement. Since $r$
can be continued to an bounded $A$-linear mapping $r(\cdot)=\langle.,x_r
\rangle$ on the $A$-dual (self-dual) Hilbert $A$-module ${\cal M}'$ of
${\cal M}$ ($x_r \in {\rm Ker}(r)^\bot \subseteq {\cal M}'$) and since the
orthogonal complement of the kernel of $r$ inside ${\cal M}'$ is a direct
summand isomorphic to $\{ Ap, \langle .,. \rangle \}$ for some projection
$p \in A$ by the structural theorem \cite[Th.~3.12]{Pa} for self-dual
Hilbert W*-modules the orthogonal complement of the kernel of $r$ with
respect to ${\cal M}$ is isomorphic to the Hilbert $A$-module $\{ I,\langle
.,. \rangle_A \}$ for some norm-closed left ideal $I \subseteq Ap$ of $A$,
where the left-strict closure of the left ideal $I$ is the w*-closed ideal
$Ap$ of $A$. Now, $r$ can be identified with the element $x_r \in Ap$, and
$x_r \in Ap$ is the left-strict limit of a net $\{ x_\beta : \beta \in J \}$
of elements of $I \cap {\cal M}$, cf.~\cite[\S 3.12]{Ped}.

\noindent
Finally, by transfinit induction one has to decompose ${\cal M}$ into a sum
of pairwise orthogonal direct summands of type ${\rm Ker}(r)^\bot$ for
bounded $A$-linear functionals $r$ on ${\cal M}$, where ${\rm Ker}(r)^\bot$
is always isomorphic to a left norm-closed ideal $I$ of $A$ with the
standard $A$-valued inner product on it.
$\: \bullet$

\medskip \noindent
We go on to investigate the image of bounded module mappings between
Hilbert W*-modules. In general, many quite non-regular things can happen as
the example below shows, but embeddings of self-dual Hilbert W*-modules
into other Hilbert W*-modules can be shown to be mappings onto direct
summands in contrast to the situation for general Hilbert C*-modules.

\begin{ex}  \label{exam2} {\rm
    Let $A$ be the set of all bounded linear operators $B(H)$ on a separable
    Hilbert space $H$ with basis $\{ e_i : i \in {\bf N}\}$. Denote by $k$
    the operator $k(e_i)=\lambda_i e_i$ for a sequence $\{ \lambda_i : i \in
    {\bf N} \}$ of non-zero positive real numbers converging to zero. Then
    the mapping
    \[
    \phi_k: A \rightarrow A \quad , \quad \phi_k: a \rightarrow a \cdot k
    \]
    is a bounded $A$-linear mapping on the left projective Hilbert
    $A$-module $A$. But the image is not a direct summand of this $A$-module
    and is not even Hilbert because direct summands of $A$ are of the form
    $Ap$ for some projection $p$ of $A$, and $1_A \cdot k$ should equal $p$.
    The image of $\phi_k$ is a subset of the set of all compact operators on
    $H$. Note, that the mapping $\phi_k$ is not injective.}
\end{ex}

\noindent
The following proposition was proved for arbitrary C*-algebras $A$, countably
generated Hilbert $A$-modules ${\cal M}$, ${\cal N}$ without self-duality
restriction and an injective bounded module mapping $\alpha:{\cal M}
\rightarrow {\cal N}$ with norm-dense range by H.~Lin \cite[Th.~2.2]{Lin:92}.
We present another variant for a similar situation in the W*-case.

\begin{prop}     \label{prop2}
    Let $A$ be a W*-algebra, ${\cal M}$ be a self-dual Hilbert $A$-module
    and $\{ {\cal N}, \langle .,. \rangle \}$ be another Hilbert $A$-module.
    Suppose, there exists an injective bounded module mapping $\alpha:
    {\cal M} \rightarrow {\cal N}$ with the range property $\alpha
    ({\cal M})^{\bot \bot} = {\cal N}$.
    Then the operator $\alpha(\alpha^*\alpha)^{-1/2}$ is a bounded
    module isomorphism of ${\cal M}$ and ${\cal N}$.
    In particular, they are isomorphic as Hilbert $A$-modules.
\end{prop}

\noindent
{\it Proof:\/} The mapping $\alpha$ possesses an adjoint bounded module
mapping $\alpha^*:{\cal N} \rightarrow {\cal M}$ because of the self-duality
of ${\cal M}$, cf.~\cite[Prop. 3.4]{Pa}.
Since $\alpha^*\alpha$ is a positive element of the C*-algebra ${\rm End}_A
({\cal M})$ of all bounded (adjointable) module mappings on the Hilbert
$A$-module ${\cal M}$ the square root of it, $(\alpha^*\alpha)^{1/2}$,
is well-defined by the series
\[
(\alpha^*\alpha)^{1/2}=
\|.\|-\lim_{n \rightarrow \infty} \|(\alpha^*\alpha)\|^{1/2}
\left(
{\rm id}_{\cal M} -
\sum_{k=1}^n \lambda_k \left( {\rm id}_{\cal M} -
\frac{(\alpha^*\alpha)}{\|(\alpha^*\alpha)\|}
\right)^k\right)
\]
with coefficients $\{ \lambda_k \}$ taken from the Taylor series at
zero of the complex-valued function $f(x)=\sqrt{1-x}$ on the interval
[0,1]. Moreover, because
 \[
\langle (\alpha^*\alpha)^{1/2}(x),(\alpha^*\alpha)^{1/2}(x) \rangle =
\langle \alpha(x),\alpha(x) \rangle
\]
and because of the injectivity of $\alpha$ the mapping
$(\alpha^*\alpha)^{1/2}$ has trivial kernel. At the contrary one can only
say that the range of $(\alpha^*\alpha)^{1/2}$ is $\tau_1$-dense in
${\cal M}$, (cf.~\cite{Frank:ZAA}). Indeed, for every $A$-linear bounded
functional $r(\cdot )=\langle .,y \rangle$ on the self-dual Hilbert
$A$-module ${\cal M}$ mapping the range of $(\alpha^*\alpha)^{1/2}$
to zero one has
\[
0 = \langle (\alpha^*\alpha)^{1/2}(x),y \rangle
= \langle x, (\alpha^*\alpha)^{1/2}(y) \rangle
\]
for every $x \in {\cal M}$. Hence, $y=0$ since $(\alpha^*\alpha)^{1/2}$
is injective and $x \in {\cal M}$ was arbitrarily chosen.

\noindent
Now, consider the mapping $\alpha(\alpha^*\alpha)^{-1/2}$ where it is defined
on ${\cal M}$. Since $(\alpha^*\alpha)^{1/2}$ has both $\tau_1$-dense range
and trivial kernel by the assumptions on $\alpha$ its inverse unbounded
module operator $(\alpha^*\alpha)^{-1/2}$ is $\tau_1$-densely defined.
One obtains
\[
\langle
\alpha(\alpha^*\alpha)^{-1/2}(x),\alpha(\alpha^*\alpha)^{-1/2}(y)
\rangle
 =  \langle x,y \rangle
\]
for every $x,y$ from the ($\tau_1$-dense) area of definition of
$(\alpha^*\alpha)^{-1/2}$.
Consequently, the operator $\alpha(\alpha^*\alpha)^{-1/2}$ continues
to a bounded isometric module operator on ${\cal M}$ by
$\tau_1$-continuity. The range of it is $\tau_1$-closed (i.e., a
self-dual direct summand of ${\cal N}$) and hence, equals ${\cal N}$ by
assumption. Finally, since the range of $(\alpha^*\alpha)^{-1/2}$ is
norm-closed and $\tau_1$-dense in $\cal M$
and since $\cal M$ is self-dual
the mapping $\alpha$ is a (non-isometric, in general)
Hilbert $A$-module isomorphism itself.
$\: \bullet$

\begin{cor}  \label{corol2}
    Let $A$ be a W*-algebra, ${\cal M}$ be a self-dual Hilbert
    $A$-module and $\{ {\cal N}, \langle .,. \rangle \}$ be another
    Hilbert $A$-module. Every injective module mapping from ${\cal M}$ into
    ${\cal N}$ is a Hilbert $A$-module isomorphism of ${\cal M}$ and
    of a direct summand of ${\cal N}$.
\end{cor}

\noindent
For our application in \S 3 we need the following partial
result:

\begin{cor}  \label{corol2a}
    Let $A$ be a W*-algebra, ${\cal M}$ and ${\cal N}$ be countably generated
    Hilbert $A$-modules and $F:{\cal M} \to {\cal N}$ be a
    Fredholm operator (see {\rm ~\cite{MF}}). Then $\Ker F$ and
    $(\Im F)^{\bot}$ are projective finitely generated $A$-submodules,
    and $\Ind F=[\Ker F] -[(\Im F)^{\bot}]$ inside $K_0(A)$.
\end{cor}

\noindent
{\it Proof:\/} We denote by $\hop$ the direct orthogonal sum of two
Hilbert $A$-modules, whereas $\bop$ denotes the direct topological
sum of two Hilbert $A$-submodules of a given Hilbert $A$-module,
where orthogonality of the two components is not required.
Let ${\cal M}={\cal M}_0 \hop {\cal M}_1$,
${\cal N}={\cal N}_0 \bop {\cal N}_1$
be the decompositions from the definition of $A$-Fredholm operator:
\[
F=\Mat {F_0}00{F_1} :
{\cal M}_0 \hop {\cal M}_1   \to   {\cal N}_0 \bop {\cal N}_1,
\]
$F_0: {\cal M}_0 \cong {\cal N}_0$, $F_1:{\cal M}_1 \to {\cal N}_1$,
${\cal M}_1$ and ${\cal N}_1$ are the projective finitely generated
modules. Let
$x=x_0+x_1$, $x_0\in {\cal M}_0$ and $x_1\in {\cal M}_1$, and
$F(x)=0$, so $0=F_0(x_0)+F_1(x_1)\in {\cal N}_0\bop{\cal N}_1$. Thus
$F_0(x_0)=0$, $F_1(x_1)=0$, so $x_0=0$ and $x\in {\cal M}_1$. Thus
$\Ker F=\Ker F_1\ss {\cal M}_1$. By Lemma~\ref{lem2} $\Ker F$ is a
projective finitely generated $A$-module and has an orthogonal
complement. So, by Corollary~\ref{corol2},
\[
F=\Matt {F_0}000{F'_1}0000 :{\cal M}_0\hop{\cal M}'_1\bop \Ker F
\to \left({\cal N}_0\bop \overline{F({\cal M}'_1)} \right) \hop (\Im
F)^\bot
\]
and $\Ind F =[\Ker F] - [(\Im F)^\bot]$.
$\bullet$

\smallskip \noindent
The following example shows that the situations may be quite different for
general Hilbert C*-modules and injective mappings between them:

\begin{ex}  \label{exam3}  {\rm
    Consider the C*-algebra $A=C([0,1])$ of all continuous functions on the
    interval [0,1] as a self-dual Hilbert $A$-module over itself equipped
    with the standard $A$-valued inner product $\langle a,b \rangle_A=ab^*$.
    The mapping $\phi: f(x) \rightarrow x \cdot f(x)$, $(x \in [0,1]$), is
    an injective bounded module mapping. Its range has trivial orthogonal
    complement, but it is not closed in norm and, consequently, not a direct
    summand of $A$. Nevertheless, the bi-orthogonal complement of the range
    of $\phi$ with respect to $A$ equals $A$.}
\end{ex}

\begin{lem}   \label{prointersect}
    Let $A$ be a W*-algebra. Let ${\cal P}$ and ${\cal Q}$ be self-dual
    Hilbert $A$-submodules of a Hilbert $A$-module ${\cal M}$. Then
    ${\cal P} \cap {\cal Q}$ is a self-dual Hilbert $A$-module and direct
    summand of ${\cal M}$. Moreover, ${\cal P}+{\cal Q}
    \subseteq {\cal M}$ is a self-dual Hilbert $A$-submodule.
\newline
    If ${\cal P}$ is projective and finitely generated then the intersection
    ${\cal P} \cap {\cal Q}$ is projective and finitely generated, too. If
    both ${\cal P}$ and ${\cal Q}$ are projective and finitely generated then
    the sum ${\cal P}+{\cal Q}$ is also.
\end{lem}

\noindent
{\em Proof.\/}
Let $p:{\cal M}={\cal P} \bop {\cal P}^\bot \rightarrow {\cal P}^\bot$ be
the canonical orthogonal projection existing by \cite[Th.~2.8]{Frank:ZAA},
(cf.~\cite{DupFil} for the projective case). Let $p_Q=p : {\cal Q} \to
{\cal P}^\bot$. Since ${\cal Q}$ is a self-dual Hilbert
$A$-module $p_Q$ admits an adjoint operator and ${\rm Ker}p_Q \subseteq
{\cal Q}$ is a direct summand by Lemma~\ref{lem2}. Consequently, it is a
self-dual Hilbert $A$-submodule of ${\cal Q} \subseteq {\cal M}$. But
${\rm Ker}p_Q = {\cal P} \cap {\cal Q}$.
To derive the second assertion one has to apply the fact again that
every self-dual Hilbert $A$-submodule is a direct
summand, cf.~\cite{Frank:ZAA}. \newline
If ${\cal P}$ is projective and finitely generated then every direct
summand of it is projective and finitely generated, what shows the
additional remarks.
$\bullet$


\section{Lefschetz numbers}        \label{sec-lef}

\noindent
Throughout this section $A$ denotes a W*-algebra. This restriction
enables us to apply the results of the previuos section being valid
only in the W*-case, in general.

\noindent
Let $U$ be a unitary operator in the projective finitely
generated Hilbert $A$-module $\cal P$. Then
  \begin{equation}\label{eq-specth}
     U=\int_{S^1}\, e^{i\vfi}\, dP(\vfi ),                 
  \end{equation}
where $P(\vfi)$ is the projection valued measure valued in the
W*-algebra of all bounded (adjointable) module operators on
$\cal P$, and the integral converges with respect to the norm.
So we have a bounded and measurable function
  \begin{equation}   \label{specfun}
    L({\cal P}, U):\, S^1 \to K_0(A),\,\vfi\mapsto [dP(\vfi )],
  \end{equation}
This function is bounded in the sense that there exists a projection
which is greater then all values with respect to the partial order in
the space of projections. Let us denote the set of such functions by
$K_0(A)_S$.

\noindent
Let us note that the Lefschetz numbers for compact group action
considered in ~\cite{TroBoch} can be thought of as evaluated (for
unitary representation) in the subspace of finitely valued (simple)
functions:
\[
{\rm Simple}(S^1,K_0(A)) \ss K_0(A)_S.
\]
Suppose, ${\cal P}=A^n$. In the case of $L({\cal P},U) \in
{\rm Simple}(S^1,K_0(A))$ associate with the integral
(\ref{eq-specth})
\[
     \int_{S^1}\, e^{i\vfi}\, dP(\vfi )
     = \sum _k\,e^{i\vfi _k}P({\cal E}_k)
\]
the following class of the cyclic homology $HC_{2l}(M(n,A))$:
\[
     \sum_k\,P({\cal E}_k) \o \dots \o P({\cal E}_k) \cdot e^{i\vfi _k}.
\]
Passing to the limit we get the following element
\[
     \tilde T U=\int_{S^1}\,e^{i\vfi} \, d(P\o \dots \o P)(\vfi)
     \in HC_{2l} (M(n,A)).
\]
Then we define
\[
     T(U)=\Tr ^n_*(\tilde T U) \,\in\,HC_{2l}(A),
\]
where $\Tr ^n_*$ is the trace in cyclic homology.

\begin{lem} \label{le-uni} {\bf (~\cite[Lemma 6.1]{TroBoch})\/}
\newline
  Let $J:\,{\cal M}=A^m \to {\cal N}=A^n$ be an isomorphism,
  $U_{\cal M}:\,{\cal M}\to {\cal M},\,U_{\cal N}:\,{\cal N}\to
  {\cal N}$ be $A$-unitary operators and
  $JU_{\cal M}=U_{\cal N}J$. Then
  \[
     T(U_{\cal M})=T(U_{\cal N}).
  \]
\end{lem}

\noindent
Similar techniques can be developed for a projective finitely
generated $A$-module ${\cal N}$ instead of $A^n$. For this
purpose we take ${\cal N}=q(A^n)$, where $q$ denotes the
orthogonal projection from $A^n$ onto its direct orthogonal
summand $\cal N$. Then we set
\[
   U\bop 1:  A^n  \cong  {\cal N} \bop  (1-q)A^n \to
           {\cal N} \bop (1-q)A^n \cong A^n,
\]
\[
   \tilde T U=\int_{S^1}\,e^{i\vfi}\,d(qPq\o\dots\o qPq)(\vfi) \, .
\]
The correctness is an immediate consequence of the Lemma \ref{le-uni}.

\smallskip \noindent
Let us consider an $A$-elliptic complex $(E,d)$ and its unitary
endomorphism $U$. The results of \S 1 (cf.~Prop.~\ref{prop2}, Lemma
\ref{corol2a}, Lemma \ref{prointersect}) and the standard Hodge theory
argument help us to prove the following lemma.

\begin{lem}   \label{le-hodraz}
  For the $A$-Fredholm operator
  \[
       F=d+d^*:\G({\cal E}_{ev}) \to \G({\cal E}_{od}),
  \]
  we have
  \[
    \Ker (F\v _{\G({\cal E}_{ev})}) \stackrel{{\rm def}}{=}
    H_{ev}({\cal E})= \bop H_{2i}({\cal E}),
  \]
  \[
    \Ker (F\v _{\G({\cal E}_{od})}) \stackrel{{\rm def}}{=}
    H_{od}({\cal E})= \bop H_{2i+1}({\cal E}),
  \]
  where $H_m({\cal E})$ is the orthogonal complement to $\Im d \ss \Ker d \ss
  \G({\cal E}_m)$ and $H_m({\cal E})$ are projective $U$-invariant
  Hilbert $A$-modules.
\end{lem}

\noindent
{\em Proof.\/} For $u_{2i}\in \G(E_{2i})$ while
$$
(d+d^*)(u_0+u_2+u_4+\dots ) =0
$$
we have
$$
du_0+d^*u_2=0,\,du_2+d^*u_4=0,\dots
$$
Together with the equality
$$
(du,d^*v)=(d^2u,v)=0
$$
one obtains
$$
du_0=0,\, du_2=0,\dots \quad ; d^*u_2=0,\,d^*u_4=0,\dots.
$$
what implies $u_{2i}\in \Ker (d+d^*)$. On the other hand for
$v_2\in\Im d,\, v_2=dv_1$ we have
$$
(v_2,u_2)=(dv_1,u_2)=(v_1,d^*u_2)=0.
$$
Thus $u_{2i}\in H_{2i}(E)$.
Conversely, let $u=u_0+u_2+\dots,\, u_{2i}\in H_{2i}(E)$, i.e.
$du_{2i}=0,\,(i=0,1,2,\dots)$, and for any $v_{2i-1}\in E_{2i-1}$ we
have
$$
(dv_{2i-1},u_{2i})=0,\quad (v_{2i-1},d^*u_{2i})=0,
$$
so $d^*u_{2i}=0$. Thus $u\in \Ker (d+d^*)$.
The invariance and projectivity follow from the proved identification
and Corollary ~\ref{corol2a}.
$\bullet$

\begin{dfn} {\rm
     We define {\it the Lefschetz number\/} $L_1$ as
     $$
       L_{1} ({\cal E},U)=
       \sum_i (-1)^i\,T(U \v H_i({\cal E})) \in K_0(A)_S.
     $$  }
\end{dfn}

\begin{dfn}  {\rm
     We define {\it the Lefschetz number\/} $L_{2l}$ as
     $$
       L_{2l} ({\cal E},U)=
       \sum_i (-1)^i\,T(U \v H_i({\cal E})) \in HC_{2l}(A).
     $$  }
\end{dfn}

\noindent
After all the following theorem is evident:

\begin{th}
  Let the Chern character $\Ch$ be defined as in
  {\rm ~\cite{ConNC,KarHC,KarHCK1}}. Then
  $$
     L_{2l}({\cal E},U)=
     \int_{S^1}({\Ch} ^0_{2l})_*(L_1({\cal E},U))(\vfi)\,d\vfi.
  $$
\end{th}

\noindent
\begin{rem}   {\rm
In situations, when the endomorphism $V$ of
the elliptic C*-complex represents as an element of a
represented there amenable group $G$ acting on the C*-complex
then the $A$-valued inner products can be chosen $G$-invariant,
what gives us the unitarity of $V$ (see ~\cite{Manu}). However,
there is another obstruction demanding new approaches which will
be shown at Example \ref{ex993} below.}
\end{rem}

\section{Obstructions in the C*-case and related topics}
\label{sec-obstr}

\noindent
The aim of this chapter is to show some obstructions arising in the
general Hilbert C*-module theory for more general C*-algebras than
W*-algebras which cause the made restriction of the investigations in
section three. The results underline the outstanding properties of
Hilbert W*-modules.
To handle the general C*-case we often need a basic construction
introduced by W.~L.~Paschke and H.~Lin. It gives a link between
the W*-case and the general C*-case.

\begin{rem}  \label{***} {\rm  (cf.\cite[Def.~1.3]{Lin:90/2},
                             \cite[\S 4]{Pa})   \newline
    Let $\{ {\cal M},\langle .,. \rangle \}$ be a left pre-Hilbert
    $A$-module over a fixed C*-algebra $A$. The algebraic tensor
    product $A^{**}\otimes {\cal M}$ becomes a left
    $A^{**}$-module defining the action of $A^{**}$ on its
    elementary tensors by the formula $ab \otimes h = a (b \otimes h)$ for
    $a,b \in A^{**}$, $h \in \cal M$. Now , setting
    \[
    \left[ \sum_i a_i \otimes h_i , \sum_j b_j \otimes g_j \right] =
    \sum_{i,j} a_i \langle h_i, g_j \rangle b_j
    \]
    on finite sums of elementary tensors one obtains a degenerate
    $A^{**}$-valued inner pre-product. Factorizing $A^{**}
    \otimes {\cal M}$ by $N= \{ z \in A^{**} \otimes {\cal M} : [z,z]=0
    \}$ one obtains a pre-Hilbert $A^{**}$-module denoted by
    ${\cal M}^{\#}$ in the sequel. It contains $\cal M$ as a $A$-submodule.
    If $\cal M$ is Hilbert then ${\cal M}^{\#}$ is Hilbert, and vice versa.
    The transfer of the self-duality is more difficult. If $\cal M$ is
    self-dual then ${\cal M}^{\#}$ is self-dual, too. But,

\noindent
\begin{prb} {\rm
    Suppose, the underlying C*-algebra $A$ is unital.
    Whether the property of ${\cal M}^{\#}$ to be self-dual implies that
    $\cal M$ was already self-dual?}
\end{prb}

\noindent
    Other standard properties like e.g.~C*-reflexivity can not be transferred.
    But every boun\-ded $A$-linear operator $T$ on $\cal M$ has a unique
    extension to a bounded $A^{**}$-linear operator on ${\cal M}^{\#}$
    preserving the operator norm, (cf.~\cite[Def.~1.3]{Lin:90/2}).
}
\end{rem}

\begin{prop} \label{prop3}
    Let $A$ be a C*-algebra, ${\cal M}$ and ${\cal N}$ be two Hilbert
    $A$-modules and $\phi : {\cal M} \rightarrow {\cal N}$ be a bounded
    $A$-linear mapping. Then the kernel ${\rm Ker}(\phi)$ of $\phi$ coincides
    with its bi-orthogonal complement inside ${\cal M}$. In general, it is
    not a direct summand.
\end{prop}

\noindent
{\it Proof: \/} Let us assume, ${\rm Ker}(\phi) \not=
{\rm Ker}(\phi)^{\bot\bot}$ with respect to the $A$-valued inner product of
${\cal M}$. Form the direct sum ${\cal L}={\cal M} \oplus {\cal N}$.
The mapping $\phi$ extends to a bounded $A$-linear mapping $\psi$ on
${\cal L}$ setting
\[
\psi(x) = \left\{ \begin{array}{r@{\quad : \quad}l}
                  \phi(x) &  x \in {\cal M} \\
                  0       &  x \in {\cal N}
                  \end{array} \right.    \, .
\]
Extend $\psi$ further to a bounded $A^{**}$-linear operator on the
correspondent Hilbert $A^{**}$-module ${\cal L}^\#$. By Lemma \ref{lem2}
the sets ${\rm Ker}(\phi)^{\#}$ and $({\rm Ker}(\phi)^{\bot\bot})^{\#}$
both are contained in the kernel ${\rm Ker}(\psi)$ of $\psi$, which is a
direct summand of ${\cal L}^{\#}$ and fulfils ${\rm Ker}(\psi)=
{\rm Ker}(\psi)^{\bot\bot}$. This contradicts the assumption.

\noindent
The second assertion follows from Example \ref{exam1}.
$\: \bullet$

\begin{cor} \label{corol3}
   Let $A$ be a C*-algebra and $\{ {\cal M}, \langle .,. \rangle \}$ be a
   Hilbert $A$-module. The kernel ${\rm Ker}(r)$ of every bounded module
   mapping $r: {\cal M} \rightarrow A$ coincides with its bi-orthogonal
   complement inside ${\cal M}$, but it is not a direct summand, in
   general.
\end{cor}

\begin{cor} \label{corol4}
   Let $A$ be a C*-algebra and $\{ {\cal M}, \langle .,. \rangle \}$ be a
   Hilbert $A$-module. Suppose, there exists a bounded module mapping
   $r: {\cal M} \rightarrow A$ with the property ${\rm Ker}(r)^\bot =
   \{0\}$. Then $r$ is the zero mapping.
\end{cor}

\begin{lem}  \label{lem5}
   Let $A$ be a C*-algebra and $\{ {\cal M}, \langle .,. \rangle \}$ be a
   (left) Hilbert $A$-module. For every bounded module mapping $r:{\cal M}
   \rightarrow A$ the subset ${\rm Ker}(r)^\bot \subseteq {\cal M}$ is a
   Hilbert $A$-submodule, and it is isomorphic as a Hilbert $A$-module to a
   norm-closed (left) ideal $D$ of $A$ equipped with the standard $A$-valued
   inner product $\langle .,. \rangle _A$.
\end{lem}

\noindent
{\it Proof: \/} By Corollary \ref{corol3} the set ${\rm Ker}(r)^\bot
\subseteq {\cal M}$ can be assumed to be non-zero, in general. Again,
form the Hilbert $A^{**}$-module ${\cal M}^\#$ and extend $r$ to a bounded
$A^{**}$-linear mapping $r'$ on it. The kernel of $r'$ is a direct summand
of ${\cal M}^\#$ isomorphic to a (left) norm-closed ideal of $A^{**}$ as a
Hilbert $A^{**}$-module by Corollary \ref{corol1} and Proposition \ref{prop1}.
Consequently, ${\rm Ker}(r) \subseteq {\rm Ker}(r') \cap {\cal M} \subseteq
{\cal M}^\#$ is isomorphic to a (left) norm-closed ideal $D$ of $A$ as a
(left) Hilbert $A$-module.
$\: \bullet$

\medskip \noindent
We want to get a structure theorem on the interrelation of
Hilbert C*-modules and their C*-dual Banach C*-modules. To obtain the full
picture define a new topology on (left) Hilbert C*-modules in analogy to the
(right) strict topology on C*-algebras $A$:

\begin{dfn} {\rm
Let $A$ be a C*-algebra and $\{ {\cal M}, \langle .,. \rangle \}$ be a (left)
Hilbert $A$-module. A norm-bounded net $\{ x_\alpha : \alpha \in I \}$ of
elements of ${\cal M}$ is {\it fundamental with respect to the right
$*$-strict topology} if and only if the net $\{ \langle y,x_\alpha
\rangle : \alpha \in I \}$ is a Cauchy net with respect to the norm topology
on $A$ for every $y \in {\cal M}$. The net $\{ x_\alpha : \alpha \in I \}$
{\it converges to an element $x \in {\cal M}$ with respect to the right
$*$-strict topology} if and only if
\[
\lim_{\alpha \in I} \| \langle y,x-x_\alpha \rangle \|_A =0
\]
for every $y \in {\cal M}$.   }
\end{dfn}

\newpage

\begin{th}  \label{th1}
   Let $A$ be a C*-algebra and $\{ {\cal M}, \langle .,. \rangle \}$ be a
   (left) Hilbert $A$-module. The following conditions are equivalent:
   \begin{enumerate}
   \renewcommand{\labelenumi}{(\roman{enumi})}
   \item ${\cal M}$ is self-dual.
   \item The unit ball of ${\cal M}$ is complete with respect to the
         right $*$-strict topology.
   \end{enumerate}
   Moreover, the linear hull of the completed with respect to the
   right $*$-strict topology unit ball of ${\cal M}$ coincides with the
   $A$-dual Banach $A$-module ${\cal M}'$ of ${\cal M}$.
\end{th}

\noindent
{\it Proof: \/} First, let us show the equivalence
(i)$\leftrightarrow$(ii).
Suppose the unit ball of ${\cal M}$ is complete with respect to the
right
$*$-strict topology. Consider an arbitrary non-trivial bounded module
mapping $r:{\cal M} \rightarrow A$ of norm one. Restrict the attention to the
non-zero Hilbert $A$-submodule ${\rm Ker}(r)^\bot \subseteq {\cal M}$
being isomorphic as a Hilbert $A$-module to a norm-closed (left) ideal $D$ of
$A$ equipped with the standard $A$-valued inner product $\langle .,.
\rangle _A$ by Lemma \ref{lem5}. By \cite[Th.~3.2]{Frank:ZAA} there exist
nets $\{x_\alpha : \alpha \in I \} \subset {\rm Ker}(r)^\bot$
bounded in norm by one such that $\tau_2-\lim_{\alpha \in I} x_\alpha =r$
inside the self-dual Hilbert $A^{**}$-module $(({\rm Ker}(r)^\bot)^\#)'$.
But, the values $r(y)$, $y \in {\rm Ker}(r)^\bot$, all belong to $A$ and,
in particular, to the set of all right multipliers of the C*-subalgebra and
two-sided ideal $B=\langle {\rm Ker}(r)^\bot,{\rm Ker}(r)^\bot \rangle$ of
$A$. Therefore, there exists a special net $\{x_\alpha : \alpha \in I \}
\subset {\rm Ker}(r)^\bot$ such that
\[
\|.\|_{{\cal M}}-\lim_{\alpha \in I} b(\langle y,x_\alpha \rangle
-r(y)) = 0
\]
for every $y \in A$, every $b \in B$, cf.~\cite[\S 3.12]{Ped}. Since the
set $\{ b y : b \in B, y \in {\rm Ker}(r)^\bot \}$ is norm-dense in
${\rm Ker}(r)^\bot$ one implication is shown. The opposit one follows
from the formula
\[
r(y) = \|.\|_A-\lim_{\alpha \in I} \langle y,x_\alpha \rangle \: ,
\: \, y \in {\cal M} \, ,
\]
defining a bounded module mapping $r: {\cal M} \rightarrow A$ for
every norm-bounded fundamental with respect to the right $*$-strict
topology net $\{ x_\alpha : \alpha \in I \} \in {\cal M}$. By the way
one has proved the conclusion that the $A$-dual Banach $A$-module
${\cal M}'$ of every Hilbert $A$-module ${\cal M}$ arises as the linear
hull of the completed with respect to the right $*$-strict topology
unit ball of ${\cal M}$.
$\: \bullet$

\begin{cor}  \label{corol5}
   Let $A$ be a C*-algebra and $D$ be a norm-closed (left) ideal of $A$.
   Then $\{ D, \langle .,. \rangle_A \}$ is self-dual if and only if there
   is a projection $p \in A$ such that $D \equiv Ap$ and $p \in D$.
\end{cor}

\noindent
{\it Proof: \/} If $D$ is self-dual then the identical embedding of $D$
into $A$ is a bounded $A$-linear mapping. It must be represented by an
element $p \in D$ with the property $dp^* = d$ for every $d \in D$.
That is, $pp^* = p \in D$ is positive and idempotent. The functional
property of the mapping $p$ gives the structure of $D$ as $D \equiv Ap$.
$\: \bullet$

\begin{th} \label{th2}
   Let $A$ be a C*-algebra and $\{ {\cal M}, \langle .,. \rangle \}$
   be a (left) Hilbert $A$-module. The following conditions are
   equivalent:
   \begin{enumerate}
   \item ${\cal M}$ is $A$-reflexive.
   \item Every norm bounded net $\{x_\alpha : \alpha \in I\}$ of
         elements of ${\cal M}$ for which all the nets
         $\{ r(x_\alpha): \alpha \in I \}$, $(r \in {\cal M}')$,
         converge with respect to $\|.\|_A$ has its limit $x$
         inside ${\cal M}$.
   \end{enumerate}
   Moreover, the linear hull of the completed with respect to this
   topology unit ball of ${\cal M}$ coincides with the $A$-bidual Banach
   $A$-module ${\cal M}''$ of ${\cal M}$.
\end{th}

\noindent
{\it Proof: \/} Suppose ${\cal M}$ is not self-dual because otherwise
one simply refers to Theorem \ref{th1}. Obviously, the linear hull of the
completion of the unit ball of ${\cal M}$ with respect to this topology is
a Banach $A$-module ${\cal N}$. Continue the $A$-valued inner product
from ${\cal M}$ to ${\cal N}$ by the rule
\[
\langle x,y \rangle = \lim_{\alpha \in I} \langle x_\alpha,y \rangle
\]
for every element $\langle .,y \rangle \in {\cal M}'$, where $y \in
{\cal M}$. Since the net converges with respect to the right
$*$-strict topology on ${\cal M}$, too, the limit $x$ can be
interpreted as an $A$-linear bounded
functional on ${\cal M}$. This lets to the definition of the value
$\langle x,x \rangle$ in the same manner. Consequently, ${\cal N}$ is
a Hilbert $A$-module containing ${\cal M}$ as a Hilbert $A$-submodule
and possessing the same $A$-dual Banach $A$-module ${\cal M}' \equiv
{\cal N}'$. (Cf.~\cite{Pa2} for similar constructions.) Moreover,
the unit ball of ${\cal N}$ is complete with respect to the new topology.
Since the $A$-valued inner product on ${\cal M}$ can be continued
to an $A$-valued inner product on ${\cal M}'' \equiv {\cal N}''$ by
\cite[Th.~2.4]{Pa2} every element of ${\cal M}''$ can be described
in this way, and ${\cal N}$ is $A$-reflexive.
$\: \bullet$


\begin{ex} \label{ex993} {\rm
Consider the C*-algebra $A=C([0,1])$ of all continuous functions
on the unit interval as a Hilbert $A$-module over itself. Let $U$ be
defined as
$$
   U(f)(t)=e^{it} f(t)   \, , \:   t \in [0,1]   \, ,
$$
a unitary operator. Take this unitary operator as the generator of
a unitary representation of the amenable abelian group $\Z$. All
complex irreducible representations of $\Z$ are one-dimensional.
If we would like to apply A.~S.~Mishchenko's theorem in this case
then we would have to have a finite spectrum for the generator $U$
of the representation what is not the case. Beside this, the only
projections inside $A$ and, therefore, the only self-adjoint idempotent
module operators on $A$ are $1_A$ and $0_A$, and there exists no spectral
decomposition of elements and no non-trivial direct $A$-module summand
inside $A$.   }
\end{ex}

\begin{rem}   {\rm
As it is known in all sufficient cases the morphism $S$ gives an
isomorphism of $HC_{2l}(A)$ and $HC_0(A)$ and we can work only
with the second group. In this situation we can define the
Lefschetz number $L_0\in HC_0(A)$ as in~\cite{TroLNM} for general
C*-algebras $A$.

\noindent
But for $K$-groups valued numbers even in the case of an action of
an e.g.~amenable group $G$ (see Example~\ref{ex993}) we need some kind of
infinitness and convergence, so we have to pass to $K_0(A)_S$. The
natural expression of this infinitness of eigenvalues is the spectral
decomposition, so we have to work with W*-algebras, at least for
$L_1$. The crucial moment is that in this situation there is no
theorem like~\cite{Mi}.

\noindent
Surely this argument is quite unexplicite and we have a chance for
refinement e.g.~for the monotone complete C*-algebras. But, the techniques
for the monotone complete case are rather complicated and the results
do only differ slightly from that of the W*-case, cf. ~\cite{Frank:MN}.  }
\end{rem}

\medskip \noindent
{\bf Acknowledgement.\/} The authors are indebted to Deutscher
Akademischer Austauschdienst for opening the opportunity of joint work
at the University of Leipzig in correspondence to the local
DAAD project ''Non-commutative geometry''.  \newline
The research of the second author was partially supported by the
Russian Foundation for Fundamental Research (Grant N 94-01-00108-a)
and the International Science Foundation (Grant no.~MGM000).


\noindent
\begin{tabular}{lp{0cm}l}
Fed.~Rep.~Germany  & &        Russia\\
Universit\"at Leipzig & & Moscow State University\\
FB Mathematik/Informatik & & Fakulty of Mechanics and Mathematics \\
Mathematisches Institut & & Chair of Higher Geometry and Topology\\
Augustusplatz 10 & & Vorob'ovy Gory \\
D-04109 Leipzig & & 119 899 Moscow \\
frank@mathematik.uni-leipzig.d400.de & & troitsky@difgeo.math.msu.su
\end{tabular}


\begin{thebibliography}{99}
    \bibitem{ConNC}  {\sc A.~Connes}, Non-commutative differential
       geometry, {\em Publ. Math. IHES\/} {\bf 62\/}(1985), 41-144.
    \bibitem{DupFil}  {\sc M.~J.~Dupr\'e, P.~A.~Fillmore},
       Triviality theorems for Hilbert modules, In:
       {\em Topics in modern operator theory\/}, 5th International
       conference on operator theory. Timisoara and Herculane
       (Romania), June 2--12, 1980,
       Basel-Boston-Stuttgart: Birkh\"auser Verlag, 1981,
       71--79.
    \bibitem{Frank:ZAA} {\sc M.~Frank}, Self-duality and
       C*-reflexivity of Hilbert C*-modules, {\it Zeitschr.~Anal.
       Anwendungen\/} {\bf 9\/}(1990), 165-176.
    \bibitem{Frank:MN} {\sc M.~Frank}, Hilbert C*-modules over
       monotone complete C*-algebras, to appear in {\it Mathematische
       Nachrichten}, 1995.
    \bibitem{Ha}  {\sc J.-F.~Havet}, Calcul fonctionnel continu
       dans les modules hilbertiens autoduaux, preprint,
       Universit\'e d'Orl\'eans, Orl\'eans, France, 1988.
    \bibitem{KarHC} {\sc M.~Karoubi}, Homologie cyclique des groupes
       et des alg\'ebres, {\em C. R. Ac. Sci. Paris, S\'erie 1,\/}
       {\bf 297\/}(1983), 381-384.
    \bibitem{KarHCK1} {\sc M.~Karoubi}, Homologie cyclique et
       $K$-th\'eorie alg\'ebrique. I, {\em C. R. Ac. Sci. Paris,
       S\'erie 1,\/} {\bf 297\/}(1983), no.~8, 447-450.
    \bibitem{Lance} {\sc E.~C.~Lance}, Hilbert C*-modules - a toolkit
      for operator algebraists, {\it Lecture Notes}, University of
      Leeds, School of Mathematics, Leeds, England, pp. 124, 1993.
    \bibitem{Lin:90/2} {\sc H.~Lin}, Bounded module maps and pure
       completely positive maps, {\it J.~Operator Theory} {\bf 26}(1991),
       121-138.
    \bibitem{Lin:92} {\sc H.~Lin}, Injective Hilbert C*-modules, {\it Pacific
      J.~Math. {\bf 154}}(1992), 131-164.
    \bibitem{Manu} {\sc V.~M.~Manu{\u{\i}}lov}, Representability of
       functionals and adjointability of operators on C*-Hilbert
       modules, preprint 1/94, Moscow State University, Dept.~Mech.~and
       Math., Seminar "Toplology and Analysis" , Moscow, Russia,
       Sept.~1994.
    \bibitem{Mi}   {\sc A.~S.~Mishchenko}, Representations of compact
       groups on Hilbert modules over C*-algebras (russ./engl.),
       {\it Trudy Mat.~Inst.~im.~V.~A.~Steklova \/}, {\bf 166\/}(1984),
       161-176 / {\it Proc.~Steklov Inst.~Math.\/} {\bf 166\/}(1986),
       179-195.
    \bibitem{MF} {\sc A.~S.~Mishchenko, A.~T.~Fomenko}, The index of
      elliptic operators over C*-algebras (russ./engl.),
      {\em Izv.~Akad.~Nauk SSSR, Ser.~Mat.,\/} {\bf 43\/}(1979),
      no.~4, 831-859 / {\em Math.~USSR - Izv.\/} {\bf 15\/}(1980),
      87-112.
    \bibitem{Pa}   {\sc W.~L.~Paschke}, Inner product modules over
      B*-algebras, {\it Trans. Amer. Math. Soc.\/} {\bf 182\/}(1973),
      443-468.
    \bibitem{Pa2}  {\sc W.~L.~Paschke}, The double $B$-dual of an
      inner product module over a C*-algebra $B$, {\it Canad.~J.~Math.}
      {\bf 26}(1974), 1272-1280.
    \bibitem{Ped} {\sc G.~K.~Pedersen}, ''C*-algebras and their
      automorphism groups'', Academic Press, London--New York--San
      Francisco, 1979.
    \bibitem{TroGlobA} {\sc E.~V.~Troitsky}, The index of equivariant
      elliptic operators over C*-algebras, {\em Annals Global
      Anal.~Geom.,\/} {\bf 5\/}(1987), no.~1, 3-22.
    \bibitem{TroLNM} {\sc E.~V.~Troitsky}, Lefschetz numbers of
      C*-complexes, {\em Springer Lect. Notes in Math.,\/}
      {\bf 1474\/}(1991), 193-206.
    \bibitem{TroBoch} {\sc E.~V.~Troitsky}, Some aspects of geometry
      of operators in Hilbert modules, preprint, Ruhr-Universit\"at
      Bochum, Fakult\"at f\"ur Mathematik, Bericht-Nr.~173, Jan.~1994.
\end{thebibliography}
\end{document}